\documentclass[12pt]{article}

\usepackage{amssymb}

\def\R{\mathbb{R}}

\def\C{\mathbb{C}}
\def\Q{\mathbb{Q}}

\newtheorem{theorem}{Theorem}
\newtheorem{definition}[theorem]{Definition}
\newtheorem{lemma}[theorem]{Lemma}
\newtheorem{proposition}[theorem]{Proposition}
\newtheorem{corollary}[theorem]{Corollary}

\newcommand{\bbox}{\ \hfill\rule[-1mm]{2mm}{3.2mm}}

\title{{\sc The inverse integrating factor \\ and the Poincar\'e map.} \thanks{The
authors are partially supported by a DGICYT grant number
MTM2005-06098-C02-02.} \\ {\normalsize Dedicated to Javier
Chavarriga.} }

\author{{\normalsize {\sc Isaac A. Garc\' \i a$^{\ (1)}$, H\'ector Giacomini$^{\ (2)}$ and Maite Grau$^{\ (1)}$}}}
\date{}

\begin{document}
\maketitle \vspace{-1.0cm}

\begin{abstract}
This work is concerned with planar real analytic differential
systems with an analytic inverse integrating factor defined in a
neighborhood of a regular orbit. We show that the inverse
integrating factor defines an ordinary differential equation for
the transition map along the orbit. When the regular orbit is a
limit cycle, we can determine its associated Poincar\'e return map
in terms of the inverse integrating factor. In particular, we show
that the multiplicity of a limit cycle coincides with the
vanishing multiplicity of an inverse integrating factor over it.
We also apply this result to study the homoclinic loop
bifurcation. We only consider homoclinic loops whose critical
point is a hyperbolic saddle and whose Poincar\'e return map is
not the identity. A local analysis of the inverse integrating
factor in a neighborhood of the saddle allows us to determine the
cyclicity of this polycycle in terms of the vanishing multiplicity
of an inverse integrating factor over it. Our result also applies
in the particular case in which the saddle of the homoclinic loop
is linearizable, that is, the case in which a bound for the
cyclicity of this graphic cannot be determined through an
algebraic method.
\end{abstract}
{\small{\noindent 2000 {\it AMS Subject Classification:} 37G15, 37G20, 34C05   }}

\section{Introduction and statement of the main results \label{sect1}}

We consider two--dimensional autonomous systems of real
differential equations of the form
\begin{equation}\label{gensys}
\dot{x} = P(x,y), \quad \dot{y} = Q(x,y),
\end{equation}
where $P(x,y)$ and $Q(x,y)$ are analytic functions defined on an
open set $\mathcal{U} \subseteq \mathbb{R}^2$. Here, the dot
denotes, as usual, derivative with respect to the independent
variable $t$. The vector field associated to system (\ref{gensys})
will be denoted ${\mathcal X} = P(x,y)\partial_x + Q(x,y)
\partial_y$ and its divergence is ${\rm div} \mathcal{X} :=
\partial P / \partial x + \partial Q / \partial y$.

\begin{definition}
An {\it inverse integrating factor} for system {\rm
(\ref{gensys})} in $\mathcal{U}$ is a non--locally null
$\mathcal{C}^1$ solution $V : \mathcal{U} \subset \mathbb{R}^2 \to
\mathbb{R}$ of the linear partial differential equation
\begin{equation}\label{bolasiif}
\mathcal{X} V \ = \ V \, {\rm div} \mathcal{X} \ .
\end{equation}
\end{definition}
If $V$ is an inverse integrating factor of system {\rm
(\ref{gensys})} then the zero--set of $V$, $V^{-1}(0) :=\{(x,y)
\mid V(x,y)=0 \}$, is composed of trajectories of {\rm
(\ref{gensys})}. In fact it is easy to see that for any point $p
\in \mathcal{U}$, if $\Phi(t;p)$ is the orbit of {\rm
(\ref{gensys})} that satisfies $\Phi(0;p)=p$, then
\begin{equation} \label{propdev}
V(\Phi(t;p)) = V(p)\exp\left(\int_0^t {\rm div} \mathcal{X} \circ
\Phi(s;p) \ ds \right) \ .
\end{equation}
Thus if $V(p) = 0$ then $V(\Phi(t;p)) = 0$ for all $t$ provided
that $V$ is defined on $\Phi(t;p)$.
\newline

We consider a certain regular orbit of system (\ref{gensys}) and
we assume the existence of an analytic inverse integrating factor
in a neighborhood of this regular orbit. We show that some
qualitative properties of the orbits in a neighborhood of the
considered regular orbit can be deduced from the known inverse
integrating factor.
\par
We consider a regular orbit $\phi(t)$ of system (\ref{gensys}) and
two transversal sections $\Sigma_1$ and $\Sigma_2$ based on it. We
are going to study the transition map of the flow of system
(\ref{gensys}) in a neighborhood of this regular orbit. This
transition map is studied by means of the {\em Poincar\'e map}
$\Pi: \Sigma_1 \to \Sigma_2$ which is defined as follows. Given a
point in $\Sigma_1$, we consider the orbit of system
(\ref{gensys}) with it as initial point and we follow this orbit
until it first intersects $\Sigma_2$. The map $\Pi$ makes
correspond to the point in $\Sigma_1$ the encountered point in
$\Sigma_2$.
\par
Let $(\varphi(s), \psi(s)) \in \mathcal{U}$, with $s \in
\mathcal{I} \subseteq \mathbb{R}$ be a parameterization of the
regular orbit $\phi(t)$ between the base points of $\Sigma_1$ and
$\Sigma_2$. In particular $s$ can be the time $t$ associated to
system (\ref{gensys}). Given a point $(x,y)$ in a sufficiently
small neighborhood of the orbit $(\varphi(s), \psi(s))$, we can
always encounter values of the {\em curvilinear coordinates} $(s,
n)$ that realize the following change of variables:
$x(s,n)=\varphi(s) - n \psi'(s)$, $y(s,n) =\psi(s) + n
\varphi'(s)$. We remark that the variable $n$ measures the
distance perpendicular to $\phi(t)$ from the point $(x,y)$ and,
therefore, $n=0$ corresponds to the considered regular orbit
$\phi(t)$. We can assume, without loss of generality, that the
transversal section $\Sigma_1$ corresponds to $\Sigma_1 \, := \,
\left\{ s=0 \right\}$ and $\Sigma_2$ to $\Sigma_2 \,  := \,
\left\{ s=L \right\}$, for a certain real number $L>0$.
\par
We perform the change to curvilinear coordinates $(x,y) \mapsto
(s,n)$ in a neighborhood of the regular orbit $n=0$ with $s \in
\mathcal{I} \, = \, [0, L]$. Then, system (\ref{gensys}) reads for
\begin{equation} \label{eqns.0} \dot{n} \, = \, N(s,n), \quad
\dot{s} \, = \, S(s,n), \end{equation} where $N(s,0) \equiv 0$
since $n=0$ is an orbit and $S(s,0) \neq 0$ for $s \in
\mathcal{I}$ because it is a regular orbit. Therefore, and in
order to study the behavior of the orbits in a neighborhood of
$n=0$, we can consider the following ordinary differential
equation:
\begin{equation}
\frac{dn}{ds} \, = \, F(s,n) \ . \label{eq2.0}
\end{equation}
We denote by $\Psi(s;n_0)$ the flow associated to the equation
(\ref{eq2.0}) with initial condition $\Psi(0;n_0)=n_0$. In these
coordinates, the Poincar\'e map $\Pi: \Sigma_1 \to \Sigma_2$
between these two transversal sections is given by
$\Pi(n_0)=\Psi(L;n_0)$.
\par
We assume the existence of an analytic inverse integrating factor
$V(x,y)$ in a neighborhood of the considered regular orbit
$\phi(t)$ of system (\ref{gensys}). In fact, when $\Sigma_1 \neq
\Sigma_2$ and no return is involved, there always exists such an
inverse integrating factor. It is clear that in a neighborhood of
a regular point there is always an inverse integrating factor.
Applying the characteristics' method, we can extend this function
following the flow of the system until we find a singular point or
a return is involved. \par The change to curvilinear coordinates
gives us an inverse integrating factor for equation (\ref{eq2.0}),
denoted by $\tilde{V}(s,n)$ and which satisfies
\begin{equation} \frac{\partial \tilde{V}}{\partial s} \, + \,
\frac{\partial \tilde{V}}{\partial n} \, F(s,n) \, = \,
\frac{\partial F}{\partial n} \, \tilde{V}(s,n). \label{eq3.0}
\end{equation} \par

The following theorem gives the relation between the inverse
integrating factor and the Poincar\'e map defined over the
considered regular orbit.
\begin{theorem}
We consider a regular orbit $\phi(t)$ of system {\rm
(\ref{gensys})} which has an inverse integrating factor $V(x,y)$
of class $\mathcal{C}^1$ defined in a neighborhood of it and we
consider the Poincar\'e map associated to the regular orbit
between two transversal sections $\Pi: \Sigma_1 \to \Sigma_2$. We
perform the change to curvilinear coordinates and we consider the
ordinary differential equation {\rm (\ref{eq2.0})} with the
inverse integrating factor $\tilde{V}(s,n)$ which is obtained from
$V(x,y)$. In these coordinates, the transversal sections can be
taken such that $\Sigma_1 \, := \, \left\{ s=0\right\}$ and
$\Sigma_2 \, := \, \left\{ s=L\right\}$, for a certain real value
$L>0$. We parameterize $\Sigma_1$ by a real parameter denoted
$\sigma$. The following identity holds.
\begin{equation}
\tilde{V} \left(L,\Pi(\sigma)\right) \, = \, \tilde{V}
\left(0,\sigma\right) \Pi'(\sigma). \label{eqvpi}
\end{equation}
\label{thvpi}
\end{theorem}
\vspace{-0.5cm} The proof of this result is given in Section
\ref{sect2}.
\par We remark that if we know an inverse integrating factor for system
(\ref{gensys}), we can construct the function $\tilde{V}
\left(s,n\right)$ and equation (\ref{eqvpi}) gives an ordinary
differential equation for $\Pi(\sigma)$, which is always of
separated variables. Thus, we can determine the expression of
$\Pi(\sigma)$ in an implicit way and up to quadratures. The first
example in Section \ref{sect3} illustrates this remark.  \par We
notice that the definition of the Poincar\'e map $\Pi$ is
geometric and Theorem \ref{thvpi} gives an ordinary differential
equation to compute an explicit expression of $\Pi$, provided that
an inverse integrating factor is known. As far as the authors
know, no such a way of describing the Poincar\'e map has been
given in any previous work.
\newline

It is well--known, due to the flow box theorem, that a regular
orbit of system (\ref{gensys}) which is not a separatrix of
singular points nor a limit cycle has an associated Poincar\'e map
conjugated with the identity. There are several results which
establish that the vanishing set of an inverse integrating factor
gives the orbits whose associated Poincar\'e map is not conjugated
with the identity. In this sense, Giacomini, Llibre and Viano
\cite{GLV}, showed that any limit cycle $\gamma \subset
\mathcal{U} \subseteq \mathbb{R}^2$ of system (\ref{gensys})
satisfies $\gamma \subset V^{-1}(0)$ provided that the inverse
integrating factor $V$ is defined in $\mathcal{U}$. \par Lately it
was shown that the zero--set of $V$ often contains the
separatrices of critical points in $\mathcal{U}$. More precisely,
in \cite{BerroneGiacomini} Berrone and Giacomini proved that, if
$p_0$ is a hyperbolic saddle--point of system (\ref{gensys}), then
any inverse integrating factor $V$ defined in a neighborhood of
$p_0$ vanishes on all four separatrices of $p_0$, provided
$V(p_0)=0$. We emphasize here that this result does not hold, in
general, for non--hyperbolic singularities (see \cite{Ga-Sh} for
example).
\newline

We are going to consider regular orbits whose Poincar\'e map is a
return map. We take profit from the result stated in Theorem
\ref{thvpi} in order to study the Poincar\'e map associated to a
limit cycle or to a homoclinic loop, in terms of the inverse
integrating factor. Although we have used curvilinear coordinates
to state Theorem \ref{thvpi}, we do not need to use them to
establish the results described in the following Theorems
\ref{th-mult-limv} and \ref{th-mult-loop0}. The only hypothesis
that we need is the existence and the expression of an analytic
inverse integrating factor in a neighborhood of the limit cycle or
the homoclinic loop. In particular, we are able to give the
cyclicity of a limit cycle or of a homoclinic loop from the
vanishing multiplicity of the inverse integrating factor over the
limit cycle or the homoclinic loop. The {\em vanishing
multiplicity} of an analytic inverse integrating factor $V(x,y)$
of system (\ref{gensys}) over a regular orbit $\phi(t)$ is defined
as follows. We recall the local change of coordinates
$x(s,n)=\varphi(s) - n \psi'(s)$, $y(s,n) =\psi(s) + n
\varphi'(s)$ defined in a neighborhood of the considered regular
orbit $n=0$. If we have the following Taylor development around
$n=0$:
\begin{equation} \label{eq-vnonula.0}
V(x(s,n), y(s,n)) \, = \, n^{m} \, v(s) \, + \, O(n^{m+1}) ,
\end{equation}
where $m$ is an integer with $m \geq 1$ and the function $v(s)$ is
not identically null, we say that $V$ has multiplicity $m$ on
$\phi(t)$. In fact, as we will see in Lemma \ref{lem-vnonula},
$v(s) \neq 0$ for any $s \in \mathcal{I}$, and thus, the vanishing
multiplicity of $V$ on $\phi(t)$ is well--defined over all its
points.
\newline

Let us consider as regular orbit a limit cycle $\gamma$ and  we
use the parameterization of $\gamma$ in curvilinear coordinates
$(s,n)$ with $s \in [0, L)$. The Poincar\'e return map $\Pi$
associated to $\gamma$ coincides with the previously defined map
in which $\Sigma := \Sigma_1 = \Sigma_2$. Since the qualitative
properties of $\Pi$ do not depend on the chosen transversal
section $\Sigma$, we take the one given by the curvilinear
coordinates and, thus, $\Pi(n_0) \, = \, \Psi(L;n_0)$. It is well
known that $\Pi$ is analytic in a neighborhood of $n_0=0$. We
recall that the periodic orbit $\gamma$ is a limit cycle if, and
only if, the Poincar\'e return map $\Pi$ is not the identity. If
$\Pi$ is the identity, we have that $\gamma$ belongs to a period
annulus. We recall the definition of multiplicity of a limit
cycle: $\gamma$ is said to be a limit cycle of {\em multiplicity}
$1$ if $\Pi'(0) \neq 1$ and $\gamma$ is said to be a limit cycle
of multiplicity $m$ with $m \geq 2$ if $\Pi(n_0) \, = \, n_0 \, +
\, \beta_m \, n_0^m \, + \, O(n_0^{m+1})$ with $\beta_m \neq 0$.
The following result states that given a limit cycle $\gamma$ of
(\ref{gensys}) with multiplicity $m$, then the vanishing
multiplicity of an analytic inverse integrating factor defined in
a neighborhood of it, must also be $m$.

\begin{theorem} \label{th-mult-limv}
Let $\gamma$ be a periodic orbit of system {\rm (\ref{gensys})}
and let $V$ be an analytic inverse integrating factor defined in a
neighborhood of $\gamma$.
\begin{itemize}
\item[{\rm (a)}] If $\gamma$ is a limit cycle of multiplicity $m$,
then $V$ has vanishing multiplicity $m$ on $\gamma$.

\item[{\rm (b)}] If $V$ has vanishing multiplicity $m$ on
$\gamma$, then $\gamma$ is a limit cycle of multiplicity $m$ or it
belongs to a continuum of periodic orbits.
\end{itemize}
\end{theorem}
The proof of this theorem is given in Section \ref{sect2}. \par
The idea to study the multiplicity of a limit cycle by means of
the vanishing multiplicity of the inverse integrating factor
already appears in the work \cite{GVL}. In that work the
considered limit cycles are semistable and the vanishing
multiplicity is defined in terms of polar coordinates, which only
apply for convex limit cycles. \par Since the Poincar\'e map of a
periodic orbit is an analytic function and the multiplicity of a
limit cycle is a natural number, we obtain the following corollary
from the previous result.
\begin{corollary} \label{cor-mult-limv}
Let $\gamma$ be a periodic orbit of system {\rm (\ref{gensys})}
and let $V$ be an inverse integrating factor of class
$\mathcal{C}^1$ defined in a neighborhood of $\gamma$. We take the
change to curvilinear coordinates $x(s,n)=\varphi(s) - n
\psi'(s)$, $y(s,n) =\psi(s) + n \varphi'(s)$ defined in a
neighborhood of $\gamma$. If we have that the leading term in the
following development around $n=0$:
\[
V(x(s,n), y(s,n)) \, = \, n^{\rho} \, v(s) \, + \, o(n^{\rho}) ,
\]
where $v(s) \not\equiv 0$ is such that either $\rho=0$ or $\rho
> 1$ and $\rho$ is not a natural number, then $\gamma$ belongs
to a continuum of periodic orbits.
\end{corollary}
We remark that this corollary applies, for instance, in the
following case: let us consider an invariant curve $f=0$ of system
(\ref{gensys}) with an oval $\gamma$ such that $\nabla f|_\gamma
\neq 0$ and let us assume that $\gamma$ corresponds to a periodic
orbit of the system and that there exists an inverse integrating
factor of the form $V=f^\rho g$, where $g$ is an nonzero function
of class $\mathcal{C}^1$ in a neighborhood of $\gamma$ and
$\rho\in \mathbb{R}$, with either $\rho=0$ or $\rho \geq 1$ in
order to have a $\mathcal{C}^1$ inverse integrating factor in a
neighborhood of $\gamma$. If we have that either $\rho=0$ or
$\rho$ is not a natural number, then we deduce that $\gamma$
belongs to a period annulus, by Corollary \ref{cor-mult-limv}. On
the other hand, applying Theorem \ref{th-mult-limv}, in the case
that $\rho=m \in \mathbb{N}$, we deduce that either $\gamma$
belongs to a period annulus or it is a limit cycle with
multiplicity $m$.
\newline

A regular orbit $\phi(t) = (x(t), y(t))$ of (\ref{gensys}) is
called a {\it homoclinic orbit} if $\phi(t) \to p_0$ as $t \to \pm
\infty$ for some singular point $p_0$. We emphasize that such kind
of orbits arise in the study of bifurcation phenomena as well as
in many applications in several sciences. A {\it homoclinic loop}
is the union $\Gamma = \phi(t) \cup \{p_0\}$. We assume that $p_0$
is a hyperbolic saddle, that is, a critical point of system
(\ref{gensys}) such that the eigenvalues of the Jacobian matrix $D
\mathcal{X}(p_0)$ are both real, different from zero and of
contrary sign. We remark that this type of graphics always has
associated (maybe only its inner or outer neighborhood) a
Poincar\'e return map $\Pi : \Sigma \to \Sigma$ with $\Sigma$ any
local transversal section through a regular point of $\Gamma$.
Moreover, in this work we only study compact homoclinic loops that
are the $\alpha$-- or $\omega$--limit set of the points in its
neighborhood, i.e., homoclinic loops with finite singular point
and a return map different from the identity. This fact implies
that $\Gamma$ is a compact set, that is, it does not have
intersection with the equator in the Poincar\'e compactification.
In fact the same results can be used for some homoclinic loops
whose saddle is in the equator of the Poincar\'e compactification,
but we always assume that the affine chart we are working with
completely contains the homoclinic loop, that is, the homoclinic
loops we take into account are compact in the considered affine
chart.
\par
Our goal in this work is to study the cyclicity of the described
homoclinic loop $\Gamma$ in terms of the vanishing multiplicity of
an inverse integrating factor. Roughly speaking, the {\it
cyclicity} of $\Gamma$ is the maximum number of limit cycles which
bifurcate from it under a smooth perturbation of (\ref{gensys}),
see \cite{Roussarie} for a precise definition. The study of the
cyclicity of $\Gamma$ has been tackled by the comparison between
the Poincar\'e return map of the unperturbed system (\ref{gensys})
and that one of the perturbed system in terms of the perturbation
parameters. The main result in this study is due to Roussarie and
appears in \cite{Roussarie}. Roussarie presents the asymptotic
expression of the Poincar\'e map associated to the homoclinic loop
$\Gamma$ which allows to characterize the cyclicity of $\Gamma$.
We recall this result in the following Section \ref{sect2}. A good
book on the subject where these and other results are stated with
proofs is \cite{Roussarie2}. The way we contribute to the study of
the cyclicity of $\Gamma$ is based on the use of inverse
integrating factors. We take profit from the result of Roussarie
to characterize the cyclicity of $\Gamma$ as we are able to relate
the Poincar\'e return map with the inverse integrating factor.
\par The Poincar\'e map $\Pi: \Sigma \to \Sigma$ associated to a homoclinic loop $\Gamma = \phi(t) \cup \{p_0\}$ is
given as the composition of the transition map along the regular
homoclinic orbit $\phi(t)$, which we denote by $R$, with the
transition map in a neighborhood of the saddle point $p_0$, which
we denote by $\Delta$. We recall that $\Sigma$ is a transversal
section based on a regular point of $\Gamma$. We have that $R$ is
an analytic diffeomorphism defined in a neighborhood of the
regular part $\phi(t)$ of $\Gamma$ and, in this context, it
coincides with the Poincar\'e map associated along the whole
$\phi(t)$ which satisfies the identity described in Theorem
\ref{thvpi}. The map $\Delta$ is the transition map defined in a
neighborhood of the critical point $p_0$. This map $\Delta$ is
characterized by the first non--vanishing saddle quantity
associated to $p_0$. We recall that the first saddle quantity is
$\alpha_1 \, = \, {\rm div\, } \mathcal{X}(p_0)$ and it classifies
the point $p_0$ between being strong (when $\alpha_1 \neq 0$) or
weak (when $\alpha_1=0$). If $p_0$ is a weak saddle point, the
saddle quantities are the obstructions for it to be analytically
orbitally linearizable. We give the definition of a saddle point
to be analytically orbitally linearizable in the following Section
\ref{sect2}.
\par In order to define the saddle quantities associated to $p_0$, we translate the saddle--point $p_0$ to the origin of
coordinates and we make a linear change of variables so that its
unstable (resp. stable) separatrix has the horizontal (resp.
vertical) direction at the origin. Let $p_0$ be a weak hyperbolic
saddle point situated at the origin of coordinates and whose
associated eigenvalues are taken to be $\pm 1$ by a rescaling of
time, if necessary. Then, it is well known, see for instance
\cite{Jibin}, the existence of an analytic near--identity change
of coordinates that brings the system into:
\begin{equation}\label{normal-3.0}
\begin{array}{lll}
\dot{x} & = & \displaystyle  x \, + \, \sum_{i=1}^{k-1} a_i \,
x^{i+1} y^i \, + \, a_k \, x^{k+1} y^k \, + \, \cdots \ ,
\vspace{0.2cm} \\ \dot{y} & = & \displaystyle -y \, - \,
\sum_{i=1}^{k-1} a_i \, x^{i} y^{i+1} \, - \, b_k \, x^k y^{k+1}
\, + \, \cdots \ , \end{array}
\end{equation}
with $a_k - b_k \neq 0$ and where the dots denote terms of higher
order. The first non--vanishing saddle quantity is defined by
$\alpha_{k+1} := a_k-b_k$. We remark that the first non--vanishing
saddle quantity can be obtained through an algebraic algorithm.
\newline

We consider a homoclinic loop $\Gamma = \phi(t) \cup \{p_0\}$
where $\phi(t)$ is a homoclinic orbit through the hyperbolic
saddle point $p_0$. We assume that there exists an analytic
inverse integrating factor $V$ defined in a neighborhood of the
homoclinic loop $\Gamma$. The following theorem establishes the
relation between the vanishing multiplicity of the inverse
integrating factor $V$ over the homoclinic orbit $\phi(t)$ and the
cyclicity of $\Gamma$. We remark that the vanishing multiplicity
of the inverse integrating factor also allows to determine the
first non--vanishing saddle quantity associated to the hyperbolic
saddle point $p_0$, in case it exists.
\begin{theorem} \label{th-mult-loop0}
Let $\Gamma$ be a compact homoclinic loop through the hyperbolic
saddle $p_0$ of system {\rm (\ref{gensys})} whose Poincar\'e
return map is not the identity. Let $V$ be an analytic inverse
integrating factor defined in a neighborhood of $\Gamma$ with
vanishing multiplicity $m$ over $\Gamma$. Then, $m \geq 1$ and the
first possible non--vanishing saddle quantity is $\alpha_m$.
Moreover,
\begin{itemize}
\item[{\rm (i)}] the cyclicity of $\Gamma$ is $2m-1$, if
$\alpha_m \, \neq \, 0$,

\item[{\rm (ii)}] the cyclicity of $\Gamma$ is $2m$, otherwise.
\end{itemize}
\end{theorem}
This theorem is proved in Section \ref{sect2}. Theorems
\ref{thvpi} and \ref{th-mult-loop0} are the main results presented
in this work. We consider the following situation: we have an
analytic system (\ref{gensys}) with a compact homoclinic loop
$\Gamma$ through the hyperbolic saddle point $p_0$ and whose
Poincar\'e return map is not the identity and we aim at
determining the cyclicity of $\Gamma$. Theorem \ref{th-mult-loop0}
gives us an algorithmic procedure to solve this problem, provided
we know an analytic inverse integrating factor $V$ of
(\ref{gensys}) in a neighborhood of $\Gamma$ with vanishing
multiplicity $m$ over it. We compute the $m^{th}$ saddle quantity
$\alpha_m$ at $p_0$ recalling that all the previous saddle
quantities vanish. If $\alpha_m \neq 0$ then the cyclicity of
$\Gamma$ is $2m-1$. Otherwise, when $\alpha_m = 0$, the cyclicity
of $\Gamma$ is $2m$. \par We notice that the determination of the
$m^{th}$ saddle quantity can be overcome through an algebraic
procedure. Thus, our result also applies in the particular case in
which the saddle of the homoclinic loop is linearizable, that is,
the case in which a bound for the cyclicity of this graphic cannot
be determined through an algebraic method.
\par We remark that the same result can be applied to a double
homoclinic graphic, that is, two homoclinic orbits which share the
same hyperbolic saddle point $p_0$. This kind of graphics has also
been treated in \cite{Maoan3}.
\newline

This paper is organized as follows. Next section contains the
proof of the main results, which are stated in this first section.
Moreover, in Section \ref{sect2} several additional results
related with the problems of studying the cyclicity of a limit
cycle or a homoclinic loop appear. These results allow the proof
of the main results and are interesting by themselves as they give
the state-of-the-art of the aforementioned problems. Much of the
results that we present are generalizations of previous ones and
can be better understood within the context of this work. Last
section contains several examples which illustrate and complement
our results.

\section{Additional results and proofs \label{sect2}}

We consider a regular orbit of system (\ref{gensys}) and we are
going to prove the result we have about it, Theorem \ref{thvpi}.
We consider the aforementioned transition map $\Pi: \Sigma_1 \to
\Sigma_2$ along the regular orbit. The key tool to study this
Poincar\'e map is to change the coordinates of system
(\ref{gensys}) to local coordinates in a neighborhood of the
regular orbit, that is, to change to the defined curvilinear
coordinates $(s,n)$. We recall that $n=0$ denotes the considered
regular orbit and $s \in \mathcal{I} \subset \mathbb{R}$ gives the
transition along the orbits since the two transversal sections can
be taken as $\Sigma_1 \, = \, \left\{ s:=0\right\}$ and $\Sigma_2
\, = \, \left\{ s:=L\right\}$, with $L$ a strictly positive real
value. These coordinates also allow to give a definition for the
Poincar\'e map in terms of the flow $\Psi(s;n_0)$ associated to
the ordinary differential equation (\ref{eq2.0}) as $\Pi(n_0) \, =
\, \Psi(L;n_0)$. We remark that the expression of the flow
$\Psi(s;n_0)$ in a neighborhood of the considered orbit $n_0=0$
can be encountered by means of recursive formulae at each order of
$n_0$. This recursive determination of the flow $\Psi(s;n_0)$
allows the study of the Poincar\'e map $\Pi(n_0)$. However, these
recursive formulae involve iterated integrals which can be, and
usually are, very difficult to be computed. The explanation of
this process and an application to the study of the multiplicity
of a limit cycle can be encountered in the work \cite{GaGiGr07}.
\par
We assume the existence of an inverse integrating factor $V(x,y)$
of class $\mathcal{C}^1$ defined in a neighborhood of the regular
orbit of system (\ref{gensys}). Then, we can construct an inverse
integrating factor for the equation (\ref{eq2.0}) which we denote
by $\tilde{V}(s,n)$ and which satisfies the partial differential
equation (\ref{eq3.0}). Moreover, the change to curvilinear
coordinates gives the definition of vanishing multiplicity of
$V(x,y)$ over the considered regular orbit, as described through
the expression (\ref{eq-vnonula.0}).
\par
The following auxiliary lemma gives us that the vanishing
multiplicity of $V$ over the regular orbit is well--defined on all
its points.

\begin{lemma} \label{lem-vnonula}
The function $v(s)$ appearing in {\rm (\ref{eq-vnonula.0})} is
different from zero for any $s \in \mathcal{I}$.
\end{lemma}
{\em Proof.} First, we recall that the relationship between
$V(x,y)$ and $\tilde{V}(s,n)$ is $\tilde{V}(s,n) =
V(x(s,n),y(s,n))/(J(s,n) S(s,n))$ where $J(s,n)$ is the Jacobian
of the change $(x,y) \mapsto (s,n)$ to curvilinear coordinates and
$S(s,n)$ comes from the time rescaling $t \mapsto s$ and is the
function defined by system (\ref{eqns.0}). We note that $J(s,n) =
\varphi'(s)^2 + \psi'(s)^2 + O(n)$ and thus $J(s,0) \neq 0$ for
any $s \in \mathcal{I}$. Moreover, it is clear that $S(s,0) \neq
0$ for all $s \in \mathcal{I}$.

The Taylor development of $\tilde{V}$ around $n=0$ is of the form
$$
\tilde{V}(s,n) \, = \, \tilde{v}(s)\, n^m + O(n^{m+1}) \ .
$$
Here, $\tilde{v}(s) = j(s) v(s)$ where $j(s) = 1/(J(s,0) S(s,0))
\neq 0$ for any $s \in \mathcal{I}$.

To end with, we are going to see that $\tilde{v}(s) \neq 0$ for
any $s \in \mathcal{I}$. We recall that, by assumption, $0 \in
\mathcal{I}$. Since $\tilde{V}(s,n)$ satisfies equation
(\ref{eq3.0}), we deduce that:
\begin{equation} \label{eq3t}
\tilde{V}\left(s,\Psi(s;n_0)\right) \, = \,
\tilde{V}\left(0,\Psi(0;n_0)\right)\, \exp \left\{ \int_{0}^{s}
\frac{\partial F}{\partial n}\left(\sigma,\Psi(\sigma;n_0)\right)
\, d \sigma \right\}.
\end{equation}
We develop this identity around $n_0=0$ and we get that
$$
\tilde{v}(s) \, n_0^m \, + \, O(n_0^{m+1})\, = \, \tilde{v}(0)\,
\exp \left\{ \int_{0}^{s} \frac{\partial F}{\partial
n}\left(\sigma,0 \right) \, d \sigma \right\}  n_0^m \, + \,
O(n_0^{m+1}) \ .
$$
Since $\tilde{v}(s) \not\equiv 0$, we can assume $\tilde{v}(0)
\neq 0$. Therefore, equating the coefficients of $n_0^m$ in both
members of the previous identity, we have the desired result.
\bbox
\newline

We have stated all the definitions and notation needed to give the
proof of Theorem \ref{thvpi}. \par {\em Proof of Theorem}
\ref{thvpi}. We have that $\Psi(s;n_0)$ satisfies \[
\frac{\partial \Psi}{\partial s} (s;n_0) \, = \,
F\left(s,\Psi(s;n_0)\right) \ \ \mbox{and} \ \ \Psi(0;n_0)\, = \,
n_0, \] as it is the flow of equation (\ref{eq2.0}) with initial
condition $n_0$. We differentiate the previous identities with
respect to $n_0$ and we have that:
\[ \frac{\partial}{\partial s} \left( \frac{\partial \Psi}{\partial n_0}(s;n_0)  \right)
 \, = \, \frac{\partial F}{\partial n}\left(s,\Psi(s;n_0)\right)
\left( \frac{\partial \Psi}{\partial n_0}(s;n_0)  \right) \ \
\mbox{and} \ \
 \frac{\partial \Psi}{\partial n_0}(0;n_0)  \, = \,
1. \] Hence, \[ \int_{0}^{s} \frac{\partial F}{\partial
n}\left(\sigma,\Psi(\sigma;n_0)\right) \, d \sigma \, = \,
\int_{0}^{s} \displaystyle \frac{ \frac{\partial}{\partial \sigma}
\left( \frac{\partial \Psi}{\partial n_0}(\sigma;n_0)  \right)
}{\frac{\partial \Psi}{\partial n_0}(\sigma;n_0)} \, d \sigma \, =
\, \left. \ln \left(\frac{\partial \Psi}{\partial
n_0}(\sigma;n_0)\right) \right|_{\sigma=0}^{\sigma=s}. \] Since
$\displaystyle  \frac{\partial \Psi}{\partial n_0}(0;n_0) \, = \,
1$, we deduce that
\[ \int_{0}^{s} \frac{\partial F}{\partial
n}\left(\sigma,\Psi(\sigma;n_0)\right) \, d \sigma \, = \, \ln
\left(\frac{\partial \Psi}{\partial n_0}(s;n_0)\right). \] Using
that the function $\tilde{V}(s,n)$ satisfies (\ref{eq3t}), we
conclude that
\[
\tilde{V}\left(s,\Psi(s;n_0)\right) \, = \,
\tilde{V}\left(0,\Psi(0;n_0)\right)\, \frac{\partial
\Psi}{\partial n_0}(s;n_0) \, .
\]
When we take $s=L$ and we recall that the Poincar\'e map is
defined as $\Pi(n_0) \, = \, \Psi(L;n_0)$, we deduce identity
(\ref{eqvpi}). \bbox
\newline

Let us consider as regular orbit a limit cycle $\gamma$ and  we
use the parameterization of the whole $\gamma$ in curvilinear
coordinates with $s \in [0, L)$. The value of $L$ can be taken the
length of the limit cycle $\gamma$ and $s$ can be the arc--length
parameter associated to it. The Poincar\'e return map $\Pi$
associated to $\gamma$ coincides with the previously defined map
in which $\Sigma := \Sigma_1 = \Sigma_2$. We remark that in
curvilinear coordinates $\Sigma_1$ corresponds to the value $s=0$
and $\Sigma_2$ to the value $s=L$ and that these two sections are
equal only in coordinates $(x,y)$. Since $\gamma$ is a periodic
orbit, the change to curvilinear coordinates and all the analytic
functions in $(x,y)$ in a neighborhood of $\gamma$ need to be
$L$--periodic in $s$. \par We use that $\tilde{V}(s,n)$ is
$L$--periodic in $s$ and we develop relation (\ref{eqvpi}) in a
neighborhood of $\sigma=0$ in order to get the following result.
We recall that the point $\sigma=0$ corresponds to the periodic
orbit $\gamma$ and that it is, therefore, a fixed point of its
Poincar\'e map.

\begin{proposition} \label{prop-r}
Let us consider the following Taylor development of the analytic
Poincar\'e map $\Pi(\sigma)$ associated to a periodic orbit
$\gamma$, around $\sigma=0$
$$
\Pi(\sigma) = \left\{ \begin{array}{lll} {\rm e}^{\beta_1} \,
\sigma + O(\sigma^2) & \mbox{with} & \beta_1 \neq 0 \ , \\ \sigma
+
\beta_{k} \sigma^k + O(\sigma^{k+1}) & \mbox{with} & \beta_k \neq 0 , \ k > 1 \ , \\
\sigma & \mbox{otherwise.} &
\end{array} \right.
$$
Let $m \geq 0$ be the vanishing multiplicity of an analytic
inverse integrating factor $V$ on $\gamma$.
\begin{itemize}
\item[{\rm (a)}] If $m=0$, then $\beta_k = 0$ for all $k \geq 1$.

\item[{\rm (b)}] If $m \geq 1$, then either $m$ is the lowest
subindex such that $\beta_m \neq 0$ or $\beta_k = 0$ for all $k
\geq 1$.
\end{itemize}
\end{proposition}
{\em Proof.} Let us recall the Taylor development of the function
$\tilde{V}\left(s,n_0\right)$ in a neighborhood of $n_0=0$:
$\tilde{V}\left(s,n_0\right)\, = \, \sum_{i\geq 0} \tilde{v}_i(s)
\, n_0^i$. Since the function $V(x,y)$ is not locally null, we
have that there exists an integer $m$ with $m \geq 0$ which is the
lowest index with $\tilde{v}_m(s) \neq 0$ for any $s \in
\mathcal{I}$ from Lemma \ref{lem-vnonula}. By definition, this
value of $m$ is the vanishing multiplicity of $V(x,y)$ on the
considered regular orbit. \par Since $V(x,y)$ changed to local
coordinates is $L$--periodic in $s$, we have that the function
$\tilde{V}(s,n)$ is also $L$-periodic in $s$. Then, we can write
$\tilde{V}\left(0,\sigma\right)\, = \, \sum_{i\geq 0} \nu_i \,
\sigma^i$ and $\tilde{V}\left(L,\Pi(\sigma)\right)\, = \,
\sum_{i\geq 0} \nu_i \, \Pi(\sigma)^i$. \par

If $\beta_1 \neq 0$, we evaluate identity (\ref{eqvpi}) in
$\sigma=0$ and we get that $\tilde{V}\left(L,0\right) \, = \, {\rm
e}^{\beta_1} \tilde{V}\left(0,0\right)$ which implies that
$\nu_0=0$ and, thus, the vanishing multiplicity of $V(x,y)$ on the
regular orbit is at least $1$. We develop identity (\ref{eqvpi})
in powers of $\sigma$ and the lowest order terms in both members
of the equality correspond to $\sigma^{m}$. The equation for their
coefficients is: $\nu_m {\rm e}^{m \beta_1} \, = \, \nu_m {\rm
e}^{\beta_1}$. Since $\beta_1 \neq 0$, we deduce that $m=1$.

Let $k >1$ be the lowest subindex such that $\beta_k \neq 0$. We
have that $\Pi(\sigma)\, = \sigma + \beta_{k} \sigma^k +
O(\sigma^{k+1})$. We subtract $\tilde{V}\left(0,\sigma\right)$
from both members of (\ref{eqvpi}) and we get the following
relation:
\begin{equation}
\sum_{i\geq 1} \nu_i \left(\Pi(\sigma)^i \, - \, \sigma^i\right)
\, = \, \left(\nu_0 \, + \, \sum_{i\geq 1} \nu_i \sigma^i\right)
\left(k \beta_k \sigma^{k-1} \, + \, O(\sigma^k) \right).
\label{eqa0-1}
\end{equation}
The left hand side of (\ref{eqa0-1}) has order at least $\sigma^k$
which implies that $\nu_0=0$, in order to have the same order in
both members. Therefore, $m >0$. We have that the lowest order
terms in both sides of (\ref{eqa0-1}) correspond to
$\sigma^{m+k-1}$ and the equation of their coefficients is: $m\,
\beta_k \, \nu_m = \, k \, \beta_k \, \nu_m$, which implies that
$k=m$. \bbox
\newline

{\it Proof of Theorem} \ref{th-mult-limv}. It is a straightforward
consequence of the previous Proposition \ref{prop-r}. We have
preferred to state this proposition as we are going to use the
same reasoning for the transition map associated to a homoclinic
orbit. \bbox
\newline

Let us now consider a homoclinic loop $\Gamma \, = \, \phi(t) \cup
\left\{ p_0 \right\}$ whose critical point $p_0$ is a hyperbolic
saddle--point. We will denote by $\lambda$ and $\mu$ the
eigenvalues associated to the Jacobian matrix $D \mathcal{X}(p_0)$
with the convention $\mu < 0 < \lambda$. We associate to $p_0$ its
{\it hyperbolicity ratio} $r = -\mu/\lambda$. We say that the
singular point $p_0$ is {\it strong} if ${\rm div}
\mathcal{X}(p_0) \neq 0$ (equivalently $r \neq 1$) and it is {\it
weak} otherwise. The hyperbolic saddle $p_0$ is called $p:q$ {\it
resonant} if $r = q/p \in \mathbb{Q}^+$ with $p$ and $q$ natural
and coprime numbers.
\par
A homoclinic loop $\Gamma$ is called {\it stable} ({\it unstable})
if all the trajectories in some inner or outer neighborhood of
$\Gamma$ approach $\Gamma$ as $t \to +\infty$ ($t \to -\infty$).
In the investigation of the stability of a homoclinic loop
$\Gamma$ of system (\ref{gensys}) through a saddle $p_0$ the
quantity
\begin{equation}\label{sigma}
\alpha_1 = {\rm div} \mathcal{X}(p_0)
\end{equation}
plays an important role. In short, it is well known, see for
instance p. 304 of \cite{Andronov}, that $\Gamma$ is stable
(unstable) if $\alpha_1 < 0$ ($\alpha_1 > 0$). For this reason
such kind of homoclinic loops $\Gamma$ are said to be {\it simple}
if $\alpha_1 \neq 0$ and {\it multiple} otherwise.
\par
The cyclicity of $\Gamma$ is linked with its stability. It is
notable to observe that, in the simple case $\alpha_1 \neq 0$, the
stability of $\Gamma$ is only determined by the nature of the
saddle point itself. Andronov {\em et al.} \cite{Andronov} proved
that if $\alpha_1 \neq  0$, the possible limit cycle that
bifurcates from $\Gamma$ after perturbing the system by a
multiparameter family has the same type of stability than
$\Gamma$, hence this limit cycle is unique for small values of the
parameters. After that, Cherkas \cite{Cherkas} showed that, if
$\alpha_1 = 0$ and the associated Poincar\'e map of the loop
$\Gamma$ is hyperbolic, then the maximum number of limit cycles
that can appear near $\Gamma$ perturbing the system in the
$\mathcal{C}^1$ class is 2. We recall that a real map is said to
be hyperbolic at a point if its derivative at the point has
modulus different from $1$. In \cite{Roussarie} Roussarie presents
a generalization of these results which determine the asymptotic
expression of the Poincar\'e map associated to the loop $\Gamma$.
\par
The Poincar\'e map $\Pi$ associated to $\Gamma$ is defined over a
transversal section $\Sigma$ whose base point is a regular point
of $\Gamma$. We parameterize the transversal section $\Sigma$ by
the real local coordinate $\sigma$. The value $\sigma = 0$
corresponds to the intersection of $\Sigma$ with $\Gamma$ and
$\sigma >0$ is the side of $\Gamma$ where $\Pi(\sigma)$ is
defined. Roussarie's result is twofold: on one hand the asymptotic
expansion of $\Pi(\sigma)$ is determined and, on the other hand,
the cyclicity of $\Gamma$ is deduced from it.

\begin{theorem}{\sc \cite{Roussarie}} \label{Teo-Rouss}
Let us consider any smooth perturbation of system {\rm
(\ref{gensys})}. Then we have:
\begin{itemize}
\item[(i)] If $r \neq 1$ (equivalently, $\alpha_1 \neq 0$), then
$\Pi(\sigma) = c \, \sigma^r (1+o(1))$ with $c > 0$ and at most
$1$ limit cycle can bifurcate from $\Gamma$.

\item[(ii)] If $\alpha_1=0$ and $\beta_1 \neq 0$, then
$\Pi(\sigma) = {\rm e}^{\beta_1} \sigma + o(\sigma)$ and at most
$2$ limit cycles can bifurcate from $\Gamma$.

\item[(iii)] If $\alpha_i= \beta_i = 0$ for $i=1,2,\ldots, k$ with
$k \geq 1$ and $\alpha_{k+1} \neq 0$, then $\Pi(\sigma) = \sigma +
\alpha_{k+1} \sigma^{k+1} \log \sigma + o(\sigma^{k+1} \log
\sigma)$ and at most $2 k+1$ limit cycles can bifurcate from
$\Gamma$.

\item[(iv)] If $\alpha_i = \beta_i = \alpha_k = 0$ for
$i=1,2,\ldots, k-1$ with $k > 1$ and $\beta_{k} \neq 0$, then
$\Pi(\sigma) = \sigma + \beta_k \sigma^k + o(\sigma^k)$ and at
most $2 k$ limit cycles can bifurcate from $\Gamma$.

\item[(v)] If $\alpha_i = \beta_i = 0$ for all $i \geq 1$, then
$\Pi(\sigma) = \sigma$ and the number of limit cycles that can
bifurcate from $\Gamma$ has no upper bound.
\end{itemize}
\end{theorem}

The values $\alpha_i$, with $i \geq 1$, are the saddle quantities
associated to $p_0$ and are evaluated by a local computation. The
values $\beta_i$, with $i \geq 1$, are called {\it separatrix
quantities} of $\Gamma$ and correspond to a global computation.
The determination of the $\alpha_i$ can be explicitly done through
the algebraic method of normal form theory near $p_0$.
Additionally, $\beta_1 = \int_{\Gamma} {\rm div} \mathcal{X} \
dt$, and a more complicated expression for $\beta_2$ can be
encountered in \cite{Maoan3} and it involves several iterated
integrals. On the contrary, as far as we know, there is no closed
form expression to get $\beta_i$ for $i \geq 3$. We remark that
our result provides a way to determine the cyclicity of the
homoclinic loop $\Gamma$ without computing any $\beta_i$.
\newline

A {\it graphic} $\bar{\Gamma} = \cup_{i=1}^k \phi_i(t) \cup \{
p_1,\ldots,p_k \}$ is formed by $k$ singular points $p_1, \ldots,
p_k$, $p_{k+1} =p_1$ and $k$ oriented regular orbits $\phi_1(t),
\ldots, \phi_k(t)$, connecting them such that $\phi_i(t)$ is an
unstable characteristic orbit of $p_i$ and a stable characteristic
orbit of $p_{i+1}$. A graphic may or may not have associated a
Poincar\'e return map. In case it has one, it is called a {\it
polycycle}. Of course, the homoclinic loops  that we consider are
polycycles with just one singular point. In
\cite{BerroneGiacomini}, the case of a system (\ref{gensys}) with
a compact polycycle $\bar{\Gamma}$ whose associated Poincar\'e map
is not the identity and whose critical points $p_1,p_2, \ldots,
p_k$ are all non--degenerate is studied and it is shown that if
there exists an inverse integrating factor $V$ defined in a
neighborhood of $\bar{\Gamma}$, then $\bar{\Gamma} \subseteq
V^{-1}(0)$. A strong generalization of this result is given in
\cite{Ga-Sh}, where it is showed that any compact polycycle
contained in $\mathcal{U}$ with non--identity Poincar\'e return
map is contained into the zero--set of $V$ under mild conditions.
More concretely, Garc\'{\i}a and Shafer give the following result.

\begin{theorem}{\sc \cite{Ga-Sh}}  \label{thGaSh}
Assume the existence of an analytic inverse integrating factor $V$
defined in a neighborhood $\mathcal{N}$ of any compact polycycle
$\bar{\Gamma}$ of system {\rm (\ref{gensys})} with associated
Poincar\'e return map different from the identity. Then,
$\bar{\Gamma} \subset V^{-1}(0)$.
\end{theorem}

In \cite{Ga-Sh}, this result is also given for an inverse
integrating factor $V$ with lower regularity than analytic and
assuming several conditions. In fact, a consequence of their
results is that if we consider a compact homoclinic loop $\Gamma$
whose Poincar\'e map is not the identity and such that there
exists an inverse integrating factor $V$ of class $\mathcal{C}^1$
defined in a neighborhood of $\Gamma$, then $\Gamma \subset
V^{-1}(0)$.\par Along this paper, we will work with inverse
integrating factors $V(x,y)$ analytic in a neighborhood of a
compact homoclinic loop $\Gamma = \phi(t) \cup \{ p_0 \}$ with
associated Poincar\'e map different from the identity. Hence,
$\Gamma \subset V^{-1}(0)$.
\par Regarding the existence problem, it is well known that the
partial differential equation (\ref{bolasiif}) has a solution in a
neighborhood of any regular point, but not necessarily elsewhere.
In \cite{CGGL} it is shown that if system (\ref{gensys}) is
analytic in a neighborhood $\mathcal{N} \subset \mathcal{U}$ of a
critical point that is either a strong focus, a nonresonant
hyperbolic node, or a Siegel hyperbolic saddle, then there exists
a unique analytic inverse integrating factor on $\mathcal{N}$ up
to a multiplicative constant. Hence, we only have at the moment
the local existence of inverse integrating factors in a
neighborhood of convenient singular points. But, we shall need to
know a priori whether there exists an inverse integrating factor
$V$ defined on a neighborhood of the whole homoclinic loop
$\Gamma$ of system (\ref{gensys}). This is a hard nonlocal problem
of existence of global solutions of the partial differential
equation (\ref{bolasiif}) for which we do not know its answer. We
describe one obstruction to the existence of an analytic inverse
integrating factor defined in a neighborhood of certain homoclinic
loops.

\begin{proposition} \label{prop-noV0}
Suppose that system {\rm (\ref{gensys})} has a homoclinic loop
$\Gamma$ through the hyperbolic saddle point $p_0$ which is not
orbitally linearizable, $p:q$ resonant and strong $(p \neq q)$.
Then, there is no analytic inverse integrating factor $V(x,y)$
defined in a neighborhood of $\Gamma$.
\end{proposition}
{\it Proof.} Let us assume that there exists an analytic inverse
integrating factor $V(x,y)$ defined in a neighborhood of $\Gamma$.
Let us consider $f_\lambda(x,y) = 0$ and $f_\mu(x,y) = 0$ the
local analytic expressions of each of the separatrices associated
to $p_0$, where the subindex denotes the corresponding eigenvalue.
As we will see in Proposition \ref{anul1} $(iii)$, we have that
$V(x,y)$ factorizes, as analytic function in a neighborhood of
$p_0$, in the form $V=f_\lambda^{1+k q} \, f_\mu^{1+k p} \, u$
with $u(p_0) \, \neq \, 0$ and an integer $k \geq 0$. Since
$V(x,y)$ is defined on the whole loop $\Gamma$, it needs to have
the same vanishing multiplicity on each separatrix at $p_0$, which
implies that $k=0$. Thus, we have that the local factorization of
$V$ in a neighborhood of $p_0$ is $V=f_\lambda \, f_\mu \, u$.
\par Let us take local coordinates in a neighborhood of $p_0$,
which we assume to be at the origin. We recall, see \cite{Arnold}
and the references therein, that given an analytic system
$\dot{x}=\lambda x + \cdots$, $\dot{y} = \mu y + \cdots$ near the
origin with $\mu/\lambda = -q/p \in \Q^{-}$ with $p$ and $q$
natural and coprime numbers ($p:q$ resonant saddle) then the
system has two analytic invariant curves passing through the
origin. Moreover, it is formally orbitally equivalent to
\begin{equation}
\dot{X} = p X \left[ 1 + \delta  (U^\ell + a U^{2 \ell})  \right]
\, , \quad \dot{Y} =-q Y \ , \label{eqf} \end{equation} with
$U=X^q Y^p$, $a \in \R$, $\ell$ an integer such that $\ell \geq 1$
and $\delta \in \{ 0, \pm 1 \}$. The normal form theory ensures
that $p_0$ is orbitally linearizable if, and only if, $\delta=0$.
We are assuming that $p_0$ is not orbitally linearizable and,
thus, $\delta \neq 0$.
\par
Let us consider $\tilde{V}(X,Y)$ the formal inverse integrating
factor of system (\ref{eqf}) constructed with the transformation
of $V(x,y)$ with the near--identity normalizing change of
variables divided by the jacobian of the change. We have that
$\tilde{V}(X,Y)=X \, Y \, \tilde{u}(X,Y)$ where $\tilde{u}$ is a
formal series such that $\tilde{u}(0,0) \neq 0$. We can assume,
without loss of generality, that $\tilde{u}(0,0)=1$. Easy
computations show that if $\tilde{u}$ is constant, we have that
$\tilde{V}(X,Y)$ cannot be an inverse integrating factor of system
(\ref{eqf}) with $\delta \neq 0$. If we have that $\tilde{u}$ is
not a constant, we develop $\tilde{u}$ as a formal series in $X$
and $Y$ and we can write $\tilde{u}=1+V_s(X,Y)+ \cdots$ where the
dots correspond to terms of order strictly greater than $s$ in $X$
and $Y$ and $V_s(X,Y)$ is a homogenous polynomial in $X$ and $Y$
of degree $s$. We consider the partial differential equation
satisfied by $\tilde{V}$:
\[ \begin{array}{l} \displaystyle  p X \left[ 1 + \delta  (U^\ell + a U^{2 \ell})  \right]
\frac{\partial \tilde{V}}{\partial X} \, - \, qY \,
\frac{\partial \tilde{V}}{\partial Y} \, =  \vspace{0.2cm} \\
\displaystyle \quad = \left( p-q \, + \, p \delta (1+\ell q)\,
U^{\ell} + a \delta p (1+2\ell q) \, U^{2 \ell} \right) \tilde{V}.
\end{array}
\]
We equate terms of the same lowest degree, we deduce that
$s=\ell(p+q)$ and that:
\[ pX \, \frac{\partial V_s}{\partial X} \, - \, qY \,
\frac{\partial V_s}{\partial Y} \, = \, \delta \ell p q X^{\ell q}
Y^{\ell p}. \] Easy computations show that the general solution of
this partial differential equation is $V_s(X,Y)=\delta \ell^2 p
q^2 U^{\ell} \ln X + G(U)$ where $G$ is an arbitrary function. We
deduce that this partial differential equation has no polynomial
solution $V_s(X,Y)$ unless $\delta=0$. \par We conclude that the
existence of such an inverse integrating factor $\tilde{V}(X,Y)$
implies that $\delta=0$ in contradiction with our hypothesis.
\bbox \newline

We observe that homoclinic loops considered in Proposition
\ref{prop-noV0} are simple, since the hyperbolic saddle point is
strong. As we have already mentioned in the introduction, it is
well--known that its cyclicity is $1$. Therefore, the obstruction
to the existence of an analytic inverse integrating factor in a
neighborhood of these homoclinic loops is not relevant in the
context of the bifurcation theory. \newline

The following example provides a realization of the thesis stated
in Proposition \ref{prop-noV0}, that is, we illustrate the
existence of homoclinic loops where no analytic inverse
integrating factor can exist in a neighborhood of it. We consider
the system
\begin{equation}
\dot{x}=-x+2 y+x^2 \ , \quad \dot{y}= 2 x-y-3 x^2+ \frac{3}{2} x y
\ , \label{ej-andron}
\end{equation}
studied in \cite{Andronov}. The origin is a strong saddle because
it has eigenvalues $\mu = -3$ and $\lambda = 1$ and, hence $r = 3
\neq 1$. Moreover, system (\ref{ej-andron}) possesses a homoclinic
loop $\Gamma$ through the origin contained in the invariant
algebraic curve $f(x,y)= x^2(1-x) -y^2 = 0$ and having inside a
focus. After a linear change of variables, we write system
(\ref{ej-andron}) with Jordan linear part as
$$
\dot{x}=-3 x - \frac{11}{8} x^2+2 x y- \frac{5}{8} y^2 \ , \quad
\dot{y}= y- \frac{7}{8} x^2+x y- \frac{1}{8} y^2 \ .
$$
A computation using the normal form method shows that the above
system is conjugated to the system
$$
\dot{x}=-3 x \ , \ \dot{y}= y -\frac{86579}{248832} x y^4 \ ,
$$
up to homogeneous degree 15 included. Therefore we conclude that
system (\ref{ej-andron}) is not formally orbitally linearizable
near the origin. Thus, applying Proposition \ref{prop-noV0}, it
exists no analytic inverse integrating factor $V(x,y)$ for system
(\ref{ej-andron}) defined in a neighborhood of $\Gamma$.
\par We have used the method described in \cite{ChaGiaGra} to
show that system (\ref{ej-andron}) has no irreducible invariant
algebraic curve except $x^2(1-x) -y^2 = 0$, nor exponential
factors. Therefore, from the Darboux theory of integrability it,
can be shown that this system has no Liouvillian first integral.
For the definitions and results related to invariant algebraic
curves, exponential factors and the Darboux theory of
integrability, see \cite{ChaGiaGra} and the references therein.
\newline

In this paper we will always assume as a hypothesis the existence
of an analytic inverse integrating factor $V$ defined on a
neighborhood of $\Gamma$. Under this condition, we now refer to
the uniqueness problem. In the last section of this work, we
present several examples of differential systems with a homoclinic
loop which satisfy all our hypothesis and in which an explicit
expression of an analytic inverse integrating factor is given.

\begin{proposition} \label{prop1-m}
Consider the analytic system {\rm (\ref{gensys})} having a compact
loop $\Gamma$ with a Poincar\'e return map different from the
identity. Then, assuming the existence of an analytic inverse
integrating factor $V(x,y)$ of {\rm (\ref{gensys})} defined in a
neighborhood of $\Gamma$, we have that $V(x,y)$ is unique up to a
multiplicative constant.
\end{proposition}
{\it Proof.} Assume that we have two different analytic inverse
integrating factors $V$ and $\bar{V}$ defined in a neighborhood
$\mathcal{N}$ of $\Gamma$. From Theorem \ref{thGaSh}, it follows
that $\Gamma \subset V^{-1}(0)$ and $\Gamma \subset
\bar{V}^{-1}(0)$. Let $m$ and $\bar{m}$ be the multiplicities of
vanishing of $V$ and $\bar{V}$ on $\Gamma$, respectively. We can
assume $m \geq \bar{m}$ and we can construct for system
(\ref{gensys}) the first integral $H(x,y) = V / \bar{V}$ which is
analytic in $\mathcal{N}$. Lemma \ref{lem-vnonula} ensures that
this quotient has no poles in $\mathcal{N}$. The existence of $H$
is in contradiction with the fact that $\Gamma$ is an $\alpha$--
or $\omega$--limit set with a return map different from the
identity. \bbox \newline

As a corollary of the proof of this proposition we have that if
there exists an analytic inverse integrating factor defined in a
neighborhood of a limit cycle, then it is unique up to a
multiplicative constant. \par An easy reasoning shows that if we
assume the existence of an analytic inverse integrating factor
$V(x,y)$ defined in a neighborhood of a compact loop $\Gamma$,
through the singular point $p_0$ and whose Poincar\'e return map
is not the identity, and if there exists a unique formal inverse
integrating factor $\bar{V}(x,y)$ in a neighborhood $D$ of $p_0$,
then $V(x,y)=\bar{V}(x,y)$ for all $(x,y) \in D$ up to a
multiplicative constant. \newline

Let $\Pi(\sigma)$ be the Poincar\'e return map associated to a
homoclinic loop $\Gamma$ and defined over a transversal section
$\Sigma$. Then, $\Pi$ is the composition $\Pi = R \circ \Delta$
where $R$ is an analytic diffeomorphism defined in a neighborhood
of the regular part $\phi(t)$ of $\Gamma$ and $\Delta$ is the
transition map defined in a neighborhood of the critical point
$p_0$. We parameterize $\Sigma$ by a real parameter $\sigma \geq
0$ and $\sigma=0$ corresponds to the base point. In our notation
$\sigma$ denotes both a point in $\Sigma$ and the corresponding
real value which parameterizes it. Theorem \ref{thvpi} relates the
inverse integrating factor and the diffeomorphism defined in
neighborhood of any regular orbit by the flow of system
(\ref{gensys}). Let us consider a homoclinic orbit as regular
orbit and we remark that Proposition \ref{prop-r} also applies in
the following sense. In case we consider a homoclinic orbit, the
limiting case in which the two sections tend to the saddle $p_0$
with $\Sigma_1$ following the unstable separatrix and $\Sigma_2$
the stable separatrix, gives that this Poincar\'e map coincides
with the aforementioned regular map $R$. We remark that, if
$\mathcal{I}=(a,b)$ covers the whole homoclinic orbit in the sense
of the flow, this limiting case corresponds to go from $\Sigma_1$,
based on a value of $s$ such that $s \to a^+$, to $\Sigma_2$ based
on a value of $s$ such that $s \to b^-$.
\par
We note that, since $V(x,y)$ is a well--defined function in a
neighborhood of the whole $\Gamma$, the function $\tilde{V}
\left(s, n\right)$ takes the same value in the limiting case,
i.e., when $s$ tends to the boundaries of $\mathcal{I}$. We define
$\tilde{V}\left(a,n\right) = \lim_{s \to a^+}
\tilde{V}\left(s,n\right)$ and $\tilde{V}\left(b,n\right) =
\lim_{s \to b^-} \tilde{V}\left(s,n\right)$ and these limits exist
from the same reasoning. Moreover, $\tilde{V}\left(a,n\right) =
\tilde{V}\left(b,n\right)$. Thus, the result stated in Proposition
\ref{prop-r} is also valid to study the regular map $R$ associated
to a homoclinic orbit defined with sections in the limiting case.
\newline

In order to relate inverse integrating factors and the Poincar\'e
return map associated to a homoclinic loop, we need to study the
local behavior of the solutions in a neighborhood of the critical
saddle point $p_0$. \par

The following proposition establishes some relationships between
the vanishing multiplicity of $V$ on the separatrices of a
hyperbolic saddle $p_0$ and the nature of $p_0$, provided that
$V(p_0) = 0$.
\begin{proposition} \label{anul1}
Let $V(x,y)$ be an analytic inverse integrating factor defined in
a neighborhood of a hyperbolic saddle $p_0$ with eigenvalues $\mu
< 0 < \lambda$ of an analytic system {\rm (\ref{gensys})}. Let us
consider $f_\lambda(x,y) = 0$ and $f_\mu(x,y) = 0$ the local
analytic expression of each of the separatrices associated to
$p_0$, where the subindex denotes the corresponding eigenvalue.
Then the next statements hold:
\begin{itemize}
\item[(i)] If $p_0$ is strong, then $V(p_0) \, = \, 0$.

\item[(ii)] If $p_0$ is nonresonant, then $V=f_\lambda \, f_\mu
\, u$ with $u(p_0) \, \neq \, 0$.

\item[(iii)] If $p_0$ is $p:q$ resonant and strong, then
$V=f_\lambda^{1+k q} \, f_\mu^{1+k p} \, u$ with $u(p_0) \, \neq
\, 0$ and an integer $k \geq 0$.

\item[(iv)] If $p_0$ is weak and $V(p_0) \, = \, 0$, then
$V=f_\lambda^m \, f_\mu^m \, u$ with $u(p_0) \, \neq \, 0$ and $m$
a natural number with $m \geq 1$.
\end{itemize}
\end{proposition}
{\it Proof.} Statement (i) is clear from the definition
(\ref{bolasiif}) of inverse integrating factor.

The other three statements are based on a result due to Seidenberg
\cite{Seid}, see also \cite{ChaGiaGra}. Since $V$ is analytic, we
can locally factorize near $p_0$ as $V = f_\lambda^{m_1}
f_\mu^{m_2} u$ with $u(p_0) \neq 0$ and $m_i$ nonnegative
integers. Then, we have ${\rm div} \mathcal{X}(p_0) = m_1 \lambda
+ m_2 \mu$. Hence, $(m_1-1) \lambda + (m_2-1) \mu = 0$ and
statements (ii)--(iv) easily follow. \bbox
\newline

Let us consider an analytic system (\ref{gensys}) where $p_0$ is a
weak hyperbolic saddle. By an affine change of coordinates and
rescaling the time, if necessary, the system can be written as
\begin{equation}\label{normal-1}
\dot{x} = x + f(x,y) \ , \ \dot{y} = -y + g(x,y) \ ,
\end{equation}
where $f$ and $g$ are analytic in a neighborhood of the origin
with lowest terms at least of second order. It is well known the
existence of a formal near--identity change of coordinates $(x,y)
\mapsto (X,Y)=(x+\cdots, y+\cdots)$ that brings system
(\ref{normal-1}) into the Poincar\'e normal form
\begin{equation}\label{normal-3}
\dot{X} = X \left[ 1+\sum_{i \geq 1} a_i (X Y)^i \right] \ , \
\dot{Y} = -Y \left[ 1+\sum_{i \geq 1} b_i (X Y)^i \right]  \ .
\end{equation}
From this expression, we see that system (\ref{normal-1}) has an
analytic first integral in a neighborhood of the saddle if and
only if the {\it saddle quantities} $\alpha_{i+1} := a_i-b_i$ are
zero for all $i \geq 1$, see \cite{Jibin} and references therein
for a review. In particular we observe that system
(\ref{normal-1}) has a 1:1 resonant saddle at the origin. In
\cite{Maoan}, it is proved that a planar dynamical system is
analytically orbitally linearizable at a resonant hyperbolic
saddle, that is whose hyperbolicity ratio $r$ is a rational
number, if and only if it has an analytic first integral in a
neighborhood of the saddle. We recall that a saddle is
analytically orbitally linearizable if there exists an analytic
near--identity change of coordinates transforming the system to a
local normal form such that $\dot{Y}/\dot{X}=-rY/X$ in a
neighborhood of the saddle.

We say that system (\ref{normal-1}) is analytically {\it orbitally
linearizable} (integrable) at the origin if there exists an
analytic near--identity change of coordinates transforming the
system to
$$
\dot{X} = X h(X,Y)  \ , \ \dot{Y} = -Y h(X,Y) \ ,
$$
with $h(0,0)=1$. In this case, $V(X,Y)=X^k Y^k h(X,Y)$ is a
1--parameter family of analytic inverse integrating factors for
any integer $k \geq 1$.

In \cite{Bruno}, Brjuno shows that any resonant hyperbolic saddle
point of an analytic system is analytically orbitally linearizable
if and only if it is formally orbitally linearizable. In
particular, this fact means that either there exists at least one
saddle quantity $\alpha_i$ different from zero or the system
becomes analytically orbitally linearizable. Moreover, we can
prove the following result.

\begin{theorem} \label{formal1}
Let us consider an analytic system {\rm (\ref{gensys})} with a
hyperbolic weak saddle $p_0$ whose separatrices are locally given
by $f_{1}(x,y) = 0$ and $f_{-1}(x,y) = 0$. Let $V = f_{1}^m \,
f_{-1}^m \, u$ be an analytic inverse integrating factor defined
in a neighborhood of $p_0$ such that $u(p_0) \neq 0$ and $m$ is a
nonnegative integer. If $p_0$ is not analytically orbitally
linearizable, then $m \geq 2$ and the first nonvanishing saddle
quantity is $\alpha_m$.
\end{theorem}
{\it Proof.} We can always assume that system (\ref{gensys}) takes
the form (\ref{normal-1}) where the hyperbolic weak saddle $p_0$
is at the origin and the expression of the separatrices takes the
form $f_1(x,y)=x+o(x,y)$, $f_{-1}(x,y)=y+o(x,y)$. As a result of
Proposition \ref{anul1}, the analytic inverse integrating factor
reads for $V(x,y) = f_{1}^m f_{-1}^m u$ with nonnegative integer
$m \geq 0$ and $u(0,0) \neq 0$.

If $m=0$, then $V(p_0) \neq 0$ and there exists an analytic first
integral defined on a neighborhood of $p_0$ for system
(\ref{gensys}). Hence, $p_0$ is analytically orbitally
linearizable and the saddle quantities of $p_0$ satisfy $\alpha_i
= 0$ for any $i \geq 1$.

If $m = 1$, then $V = f_{1} f_{-1} u $ with $u(p_0) \neq 0$. Then,
by statement (iii) of Theorem 5.10 of \cite{CMR}, we conclude that
$p_0$ is also analytically orbitally linearizable and the saddle
quantities of $p_0$ satisfy $\alpha_i = 0$ for any $i \geq 1$.

If $m \geq 2$, we are going to prove that the saddle quantities at
the origin satisfy $\alpha_i = 0$ for $i=1, 2, \ldots, m-1$ and
that $\alpha_m \neq 0$. We do a formal near--identity change of
coordinates $(x,y) \mapsto (X,Y)$ transforming system
(\ref{normal-1}) into the Poincar\'e normal form (\ref{normal-3}).
In these normalizing coordinates, system (\ref{normal-3}) has the
formal inverse integrating factor
\begin{equation}\label{normal-4}
\bar{V}(X,Y) = X Y \sum_{i \geq 1} \alpha_{i+1} (X Y)^i \ .
\end{equation}
We observe that functions which define system (\ref{normal-3})
give  $\bar{V}(X,Y) \, = \, \dot{X}/X\, - \, \dot{Y}/Y$. We remark
that $\bar{V}(X,Y)$ is a univariate function of the variable $X
Y$. Since there is at least one nonzero saddle quantity, system
(\ref{normal-3}) has no formal first integral near the origin. In
particular, (\ref{normal-4}) is the unique (up to multiplicative
constants) formal inverse integrating factor of system
(\ref{normal-3}). We note that performing the above formal
near--identity change of coordinates to the inverse integrating
factor $V(x,y)$ we get a formal inverse integrating factor of
system (\ref{normal-3}) of the form $X^m Y^m \bar{u}(X,Y)$, with
$\bar{u}(0,0) \neq 0$. Therefore, from uniqueness, we must have
that it coincides with the expression (\ref{normal-4}). Thus, we
have $\alpha_i = 0$ for $i=1, 2, \ldots, m-1$ and that $\alpha_m
\neq 0$. \bbox
\newline

We have analyzed the behavior of the flow near a hyperbolic saddle
point $p_0$ when an analytic inverse integrating factor $V$ is
known. In particular, we shall see that the transition map
$\Delta$ near $p_0$ can be studied by the vanishing multiplicity
of $V$ in the separatrices of $p_0$. Let us consider two
transversal sections $\Sigma_2$ and $\Sigma_1$ in a sufficiently
small neighborhood of $p_0$, where $\Sigma_2$ (resp. $\Sigma_1$)
is based on a point over the stable (resp. unstable) separatrix of
$p_0$. We parameterize $\Sigma_2$ by a real parameter $\sigma \geq
0$ and $\sigma=0$ corresponds to the base point. We recall that
the {\it transition map} $\Delta : \Sigma_2 \to \Sigma_1$ is
defined as $\Delta(\sigma) = \Phi(\tau(\sigma); \sigma) \cap
\Sigma_1$ where $\Phi(t; \sigma)$ is the flow associated to system
(\ref{gensys}) with initial condition the point $\sigma \in
\Sigma_2$ and $\tau(\sigma)$ is the minimal positive time such
that the trajectory through $\sigma$ intersects $\Sigma_1$. The
explicit asymptotic expansion of $\Delta(\sigma)$ was given by
Dulac \cite{Dulac} in terms of the hyperbolicity ratio $r$ and
saddle quantities $\alpha_i$ of $p_0$ as follows:
$$
\Delta(\sigma) = \left\{ \begin{array}{lll} c \, \sigma^r (1+o(1))
\ , \ c > 0 & \mbox{if} & r \neq 1 \ , \\ \sigma + \alpha_{k}
\sigma^k \ln \sigma + o(\sigma^k \ln \sigma) & \mbox{if} &
\alpha_1= \cdots
= \alpha_{k-1} = 0 , \ \alpha_{k} \neq 0 \ , \\
\sigma & \mbox{if} & \alpha_i = 0 , \ \mbox{for all} \ i \geq 1 .
\end{array} \right.
$$

Since the Poincar\'e return map $\Pi(\sigma)$ associated to the
compact homoclinic loop $\Gamma$ is the composition $\Pi = R \circ
\Delta$ and is different from the identity, we only have the four
possibilities (i)--(iv) described in Theorem \ref{Teo-Rouss}.
\newline

In the following theorem, the cyclicity of $\Gamma$ denotes the
maximum number of limit cycles that bifurcate from $\Gamma$ under
smooth perturbations of (\ref{gensys}).

\begin{theorem} \label{th-mult-loop}
Let $\Gamma$ be a compact homoclinic loop through the hyperbolic
saddle $p_0$ of system {\rm (\ref{gensys})} whose Poincar\'e
return map is not the identity. Let $V$ be an analytic inverse
integrating factor defined in a neighborhood of $\Gamma$ with
vanishing multiplicity $m$ over $\Gamma$. Then the following
statements hold:
\begin{itemize}
\item[{\rm (a)}] $m \geq 1$.

\item[{\rm (b)}] If $p_0$ is strong, then $m=1$ and the cyclicity
of $\Gamma$ is $1$.

\item[{\rm (c)}] If $p_0$ is weak, then:

\begin{itemize}
\item[{\rm (c.1)}] If $p_0$ is not analytically orbitally
linearizable, then $m \geq 2$, $\alpha_i= \beta_i = 0$ for
$i=1,2,\ldots, m-1$ and $\alpha_m \neq 0$. Moreover, the cyclicity
of $\Gamma$ is $2m-1$.

\item[{\rm (c.2)}] If $p_0$ is analytically orbitally
linearizable, then $\beta_1= \beta_2 = \cdots = \beta_{m-1}=0$ and
$\beta_m \neq 0$. Moreover, the cyclicity of $\Gamma$ is $2m$.
\end{itemize}
\end{itemize}
\end{theorem}

{\em Proof.} Statement (a) is a straight consequence of the fact
that $\Gamma \subset V^{-1}(0)$, see Theorem \ref{thGaSh}.

Statement (b) follows taking into account statements (ii) and
(iii) of Proposition \ref{anul1} where we recall that the
vanishing multiplicity of $V$ in each of the separatrices of $p_0$
needs to be the same since $\Gamma$ is a loop. In particular, we
observe that the value of $k$ in statement (iii) of Proposition
\ref{anul1} must be zero by the same argument.

We assume that $p_0$ is weak. We may have that all its associated
saddle quantities are zero or that there is at least one saddle
quantity different from zero. If this last case applies, the first
non--vanishing saddle quantity is $\alpha_m$ as a result of
Theorem \ref{formal1}. Therefore, we compute the $m^{th}$ saddle
quantity associated to $p_0$, knowing than the previous saddle
quantities need to be zero and we can determine if $p_0$ is
analytically orbitally linearizable (if $\alpha_m=0$) or not (if
$\alpha_m\neq 0$). The case (c.1) in the theorem corresponds to
$\alpha_m\neq 0$ and the case (c.2) to $\alpha_m=0$.

The case (c.1) is a consequence of Theorem \ref{formal1} and part
(b) of Proposition \ref{prop-r}. More concretely, since $p_0$ is
not analytically orbitally linearizable, then using Theorem
\ref{formal1} we have $m \geq 2$, $\alpha_i = 0$ for
$i=1,2,\ldots, m-1$ and $\alpha_m \neq 0$. In addition, by part
(b) of Proposition \ref{prop-r}, we get $\beta_i = 0$ for
$i=1,2,\ldots, m-1$. We remark that, in this case, either $\beta_m
\neq 0$ or $\beta_k = 0$ for all $k \geq 1$. In any case, from
statement (iii) of Theorem \ref{Teo-Rouss}, the cyclicity of
$\Gamma$ is $2 m-1$.

Finally, the proof of (c.2) works as follows. Since $p_0$ is
analytically orbitally linearizable, then $\alpha_k = 0$ for all
$k \geq 1$. By hypothesis, the Poincar\'e map is not the identity,
so $\beta_k = 0$ for all $k \geq 1$ is not possible. Hence, by
part (b) of Proposition \ref{prop-r}, we get that $\beta_i = 0$
for $i=1,2,\ldots, m-1$ and $\beta_m \neq 0$. Thus, using
statement (iv) of Theorem \ref{Teo-Rouss}, the cyclicity of
$\Gamma$ is $2 m$. \bbox \newline

{\em Proof of Theorem} \ref{th-mult-loop0}. The thesis of Theorem
\ref{th-mult-loop0} is a corollary of Theorem \ref{th-mult-loop}.
\bbox

\section{Examples \label{sect3}}

This section contains several examples which illustrate and
complete our results. The first example shows how to give an
implicit expression of the Poincar\'e map associated to a regular
orbit via identity (\ref{eqvpi}). The second example is given to
show the existence of analytic planar differential systems with a
compact homoclinic loop $\Gamma$ whose cyclicity can be determined
by the vanishing multiplicity of an inverse integrating factor
defined in a neighborhood of it. The third example targets to show
that there is no upper bound for the number of limit cycles which
bifurcate from a compact homoclinic loop whose Poincar\'e return
map is the identity. The system considered in this example has a
numerable family of inverse integrating factors such that for any
natural number $n$ there exists an analytic inverse integrating
factor whose vanishing multiplicity on the compact homoclinic loop
is $n$. \newline

{\bf Example 1.} The following system is studied in \cite{CGG}
where a complete description of its phase portrait in terms of the
parameters is given using the inverse integrating factor as a key
tool. The system
\begin{equation} \label{eqnew} \begin{array}{lll}
\dot{x} & = & \displaystyle \lambda x - y + \lambda m_1 x^3 +
(m_2-m_1+m_1m_2)x^2y+\lambda m_1 m_2 x y^2 + m_2 y^3,
\vspace{0.2cm} \\ \dot{y} & = & \displaystyle x + \lambda y - x^3
+ \lambda m_1 x^2 y + (m_1m_2-m_1-1) x y^2 + \lambda m_1 m_2 y^3,
\end{array} \end{equation}
where $\lambda$, $m_1$ and $m_2$ are arbitrary real parameters,
has the following inverse integrating factor
\begin{equation} \label{ifinew} V(x,y) \, = \, (x^2 + y^2) \, (1+
m_1 x^2 + m_1 m_2 y^2). \end{equation} In \cite{CGG} it is shown
that if $\lambda \neq 0$, $m_1 <0$, $m_1 \neq -1$ and $m_2>0$,
then the ellipse defined by $1+ m_1 x^2 + m_1 m_2 y^2 =0$, and
which we denote by $\gamma$, is a hyperbolic limit cycle of system
(\ref{eqnew}). Moreover, with the described values of the
parameters, we have that the origin of coordinates is a strong
focus whose boundary of the focal region is $\gamma$. Moreover,
since the vanishing set of $V$ is only the focus at the origin and
the limit cycle $\gamma$, we deduce by Theorem \ref{thGaSh}, that
the region, in the outside of $\gamma$ in which the orbits spiral
towards or backwards to it, is unbounded.
\par We remark that we reencounter that this limit cycle is
hyperbolic using that the vanishing multiplicity of the inverse
integrating factor over it is $1$. In this example we are not
concerned with the multiplicity of the limit cycle but with the
Poincar\'e return map associated to it. We are going to use the
ordinary differential equation stated in (\ref{eqvpi}) to give an
implicit expression for the Poincar\'e map associated to $\gamma$
in system (\ref{eqnew}). We can parameterize the ellipse by
$(\sqrt{m_2}\, \cos(s), \ \sin(s))/\sqrt{-m_1 m_2}$ with $s \in
[0,2 \pi)$ and we perform the corresponding change to curvilinear
coordinates: $(x,y) \mapsto (s,n)$ with $x=(\sqrt{m_2}-n)\,
\cos(s)/\sqrt{-m_1 m_2}$ and $y=(1-\sqrt{m_2} \, n) \,
\sin(s)/\sqrt{-m_1 m_2}$. We obtain an ordinary differential
equation of the form (\ref{eq2.0}) which describes the behavior of
the solutions of system (\ref{eqnew}) in a neighborhood of the
limit cycle (for $n=0$). The corresponding inverse integrating
factor $\tilde{V}(s,n)$ is a $2 \pi$--periodic function in $s$ and
it satisfies that: \[ \tilde{V}(0,n) \, = \, \frac{m_1\, n\,
(n-\sqrt{m_2})\, (n-2\sqrt{m_2})}{n^2-2n\sqrt{m_2} + m_2 + m_1
m_2}.
\] The ordinary differential equation (\ref{eqvpi}) for this system
reads for $ \displaystyle \tilde{V}\left(2 \pi, \Pi(\sigma)
\right) \, = \, \tilde{V}\left(0, \sigma \right) \Pi'(\sigma)$,
where $\sigma$ can be taken as a real parameter of the section
$\Sigma\, := \, \left\{ s=0 \right\}$ and such that $\sigma=0$
corresponds to the limit cycle. We observe that the point
$\sigma=\sqrt{m_2}$ corresponds to the focus point at the origin
of system (\ref{eqnew}). This differential equation can be written
in the following Pfaffian form:
\[ \frac{d \sigma}{\tilde{V}\left(0, \sigma \right)} \, = \,
\frac{d \Pi}{\tilde{V}\left(2 \pi, \Pi \right)}, \] whose
integration gives: \begin{equation} \label{eqpnew} \frac{\left(
\Pi(\sigma) \, (\Pi(\sigma)-2 \sqrt{m_2})
\right)^{(1+m_1)/2}}{(\Pi(\sigma)-\sqrt{m_2})^{m_1}} \, = \, k_0
\, \frac{\left( \sigma \, (\sigma -2 \sqrt{m_2})
\right)^{(1+m_1)/2}}{(\sigma-\sqrt{m_2})^{m_1}} , \end{equation}
where $k_0$ is an arbitrary constant of integration. We remark
that if $k_0=1$, then the function $\Pi(\sigma) \, := \, \sigma$
satisfies the implicit identity (\ref{eqpnew}). \par We observe
that the integration of (\ref{eqvpi}) always gives rise to a
Pfaffian form of separated variables. Moreover, in the case that a
return is involved we have that this Pfaffian form is symmetric in
$\Pi$ and $\sigma$, because the function $\tilde{V}(s,n)$ needs to
be $L$--periodic in $s$.
\par We are interested in the Poincar\'e map associated to the
limit cycle $\gamma$. Since it is a hyperbolic limit cycle, we
have that $\Pi(\sigma) \, = \, e^{\beta_1} \sigma \, + \,
\mathcal{O}(\sigma^2)$ with $\beta_1 \neq 0$, as described in
Proposition \ref{prop-r}. Using this expression and identity
(\ref{eqpnew}), we deduce that $k_0=\exp \left( \beta_1 \,
(m_1+1)/2\right)$. Moreover, it can be shown that $\beta_1 \, = \,
-2 \lambda T$, where $T$ is the minimal positive period of
$\gamma$ and takes the value $T=2\pi m_1/(1+m_1)$ if $m_1<-1$ and
$T=-2\pi m_1/(1+m_1)$ if $-1<m_1<0$. Therefore, $k_0=e^{-2\lambda
m_1 \pi}$ if $m_1<-1$ and $k_0=e^{2\lambda m_1 \pi}$ if
$-1<m_1<0$. We note that the inverse integrating factor determines
the multiplicity of the limit cycle $\gamma$ but not its
stability. \newline

{\bf Example 2.} We fix an integer number $m$ with $m \geq 1$ and
we consider the algebraic curve $f=0$ with
$f(x,y)=y^2-(1-x)^2(1+x)$. We have that $f=0$ has an oval in the
range $-1 \leq x \leq 1$ with a double point in $(1,0)$. We denote
by $\Gamma$ this oval. Let us consider the planar differential
system:
\begin{equation} \label{eq4.1}
\begin{array}{lll}
\displaystyle \dot{x} & = & \displaystyle - \left(
\left[(1-m)g(x,y) + f(x,y)^{m-1} \right] \frac{\partial
f}{\partial y} + f(x,y) \frac{\partial g}{\partial y} \right)
(x^2+y^2) q(x,y) \, - \, \vspace{0.2cm} \\ & & \displaystyle
\qquad f(x,y)^m \left(2(x+y)q(x,y) + (x^2+y^2)
\frac{\partial q}{\partial y}\right), \vspace{0.3cm} \\
\displaystyle \dot{y} & = & \displaystyle  \left(
\left[(1-m)g(x,y) + f(x,y)^{m-1} \right] \frac{\partial
f}{\partial x} + f(x,y) \frac{\partial g}{\partial x} \right)
(x^2+y^2) q(x,y) \, + \, \vspace{0.2cm} \\ & & \displaystyle
\qquad f(x,y)^m \left(2(x-y)q(x,y) + (x^2+y^2) \frac{\partial
q}{\partial x}\right),
\end{array}
\end{equation}
where $g(x,y)$ and $q(x,y)$ are polynomials such that the
algebraic curve $g=0$ does not intersect the oval $\Gamma$ and the
algebraic curve $q=0$ does not contain any point in the closed
region bounded by $\Gamma$. This system has $\Gamma$ as homoclinic
loop where the critical point $(1,0)$ is a hyperbolic weak saddle.
Moreover, the origin of this system is a strong focus. It can be
shown that the function $V(x,y)=(x^2+y^2) f(x,y)^m q(x,y)$ is an
inverse integrating factor, which is analytic in the whole affine
plane. Since the origin is a strong focus and $V$ is analytic in
$\mathbb{R}^2$, we have by Theorem \ref{thGaSh}, that the boundary
of the focal region must be contained in the zero--set of $V$. We
conclude that the boundary of this focal region needs to be
$\Gamma$ and, therefore, we have that $\Gamma$ is a compact
homoclinic loop through a hyperbolic saddle and whose associated
Poincar\'e return map is not the identity. Moreover, we deduce
that the vanishing multiplicity of $V$ over $\Gamma$ is $m$. The
following function \[ H(x,y) \, = \, (x+iy)^{1-i} (x-iy)^{1+i}
f(x,y) e^{g(x,y)/f(x,y)^{m-1}} q(x,y) \] is a first integral of
system (\ref{eq4.1}). In the work \cite{LliPan04} it is shown that
this form of a first integral implies that the inverse integrating
factor is polynomial and, thus, well--defined over all the real
plane. When $m=1$, $H(x,y)$ provides an analytic first integral
defined in a neighborhood of the saddle point $(1,0)$ and we
deduce that the saddle quantities associated to the critical point
$(1,0)$ are all zero. Thus, using Theorem \ref{th-mult-loop}, we
conclude that when $m=1$ the cyclicity of $\Gamma$ in system
(\ref{eq4.1}) is $2$. \par When $m>1$, using Theorem
\ref{th-mult-loop}, we deduce that the cyclicity of $\Gamma$ in
system (\ref{eq4.1}) is $2m-1$ if $\alpha_m \neq 0$ and $2m$ if
$\alpha_m=0$.
\newline

{\bf Example 3.}  Let us consider the following planar
differential system
\begin{equation} \dot{x} \, = \, -2 y, \quad \dot{y} \, = \,
-2x+3x^2,\label{eq4.2} \end{equation} which has a homoclinic loop
$\Gamma$ contained in the invariant algebraic curve $f=0$ with
$f(x,y):=y^2-x^2+x^3$. Since the origin $(0,0)$ is the saddle
contained in $\Gamma$ and it is a hyperbolic saddle, we have that
$\Gamma$ has associated a Poincar\'e return map, which is the
identity as the system is Hamiltonian with $H=f$. In addition to
the hyperbolic saddle at the origin, the system possesses a
singular point of center type in the point with coordinates $(2/3,
0)$ and no other critical point in the affine plane. We note that
any function of the form $V=f(x,y)^n$, with $n$ a natural number,
provides an analytic inverse integrating factor for the system in
the whole plane $\mathbb{R}^2$. Thus, we cannot define the
vanishing multiplicity of an inverse integrating factor on
$\Gamma$. \par We observe that for any natural number $n$, there
exist perturbations of system (\ref{eq4.2}) with at least $n$
limit cycles which bifurcate from the considered homoclinic loop
$\Gamma$, that is, there is no finite upper bound for the
cyclicity of $\Gamma$. In the following paragraph, we illustrate
this fact by exhibiting a suitable perturbation. We remark that
the fact that the vanishing multiplicity of an inverse integrating
factor on $\Gamma$ is not defined is coherent with the
nonexistence of an upper bound for the cyclicity of $\Gamma$.
\par For instance, let us fix a natural value $n$ and let
$\varepsilon$ be a nonzero real number with $|\varepsilon|$ small
enough. We take $a_i \in \mathbb{R}$, $i=1,2,\ldots,n$, such that
$a_i \neq a_j$ if $i \neq j$ and with $0<a_i\varepsilon<4/27$. For
each $i$ we have that the algebraic curve $f_i=0$ with
$f_i:=f+a_i\varepsilon$ has an oval with the point $(2/3,0)$ in
its inner region and the point $(0,0)$ in its outer region.
Moreover, if $a_i \, \varepsilon \, < \, a_j \, \varepsilon$, the
oval defined by $f_i=0$ contains the oval $f_j=0$ in its inner
region. The following system:
\[ \dot{x} \, = \, -2y, \quad \dot{y} \, = \, -2x+3x^2+\varepsilon
y \prod_{i=1}^{n}(f+a_i\varepsilon) \] is a perturbation of system
(\ref{eq4.2}) and exhibits each one of the ovals defined by the
curves $f_i=0$ as periodic orbits. Easy computations show that, if
$|\varepsilon|>0$ is small enough, the only critical points of the
perturbed system are $(0,0)$, which is a hyperbolic saddle point,
and $(2/3,0)$ which is a strong focus. It is also easy to check
that each oval described by $f_i=0$ is a hyperbolic limit cycle of
the perturbed system which bifurcates from $\Gamma$. We have, in
this way, that the considered perturbed system has at least $n$
limit cycles which bifurcate from $\Gamma$.
\newline

{\bf Acknowledgements.} We would like to thank Prof. Crist\'obal
Garc\'{\i}a from Universidad de Huelva (Spain) for his
commentaries on normal form theory and the corresponding
computations done for system (\ref{ej-andron}).
\newline

\vspace{0.5cm}

{\bf Addresses and e-mails:} \\
$^{\ (1)}$ Departament de Matem\`atica. Universitat de Lleida.
\\ Avda. Jaume II, 69. 25001 Lleida, SPAIN.
\\ {\rm E--mails:} {\tt garcia@matematica.udl.cat}, {\tt
mtgrau@matematica.udl.cat} \vspace{0.2cm}
\\
$^{\ (2)}$ Laboratoire de Math\'ematiques et Physique Th\'eorique.
C.N.R.S. UMR 6083. \\ Facult\'e des Sciences et Techniques.
Universit\'e de Tours. \\ Parc de Grandmont 37200 Tours, FRANCE.
\\ E-mail: {\tt giacomini@phys.univ-tours.fr}

\end{document}